\documentclass[11pt]{amsart}

\usepackage{amsmath,amscd,amssymb,latexsym}

\newtheorem{thm}{Theorem}[section]

\newcommand{\sm}{\left(\smallmatrix}
\newcommand{\esm}{\endsmallmatrix\right)}
\newcommand{\mat}{\begin{pmatrix}}
\newcommand{\emat}{\end{pmatrix}}

\begin{document}

\title[Exact formulas for traces of singular moduli]{Exact formulas for traces of singular moduli of higher level
 modular functions}

\author[D. Choi, D. Jeon, S.-Y. Kang \and C. H. Kim]{Dohoon Choi, Daeyeol Jeon, Soon-Yi Kang and Chang Heon Kim}

\address{Dohoon Choi, School of Mathematics, Korea Institute for Advanced Study, Seoul
130-722, Korea} \email{choija@kias.re.kr}
\address{Daeyeol Jeon, Department of
Mathematics Education, Kongju National University, Kongju 314-701,
Chungnam, Korea} \email{dyjeon@kongju.ac.kr}
\address{Soon-Yi Kang, School of Mathematics, Korea Institute for Advanced Study, Seoul
130-722, Korea} \email{sykang@kias.re.kr}
\address{Chang Heon Kim, Department of Mathematics, Seoul Women's University,
 Seoul, 139-774, Korea} \email{chkim@swu.ac.kr}

\begin{abstract}
Zagier proved that the traces of singular values of the classical
$j$-invariant are the Fourier coefficients of a weight
$\frac{3}{2}$ modular form and Duke provided a new proof of the
result by establishing an exact formula for the traces using
Niebur's work on a certain class of non-holomorphic modular forms.
In this short note, by utilizing Niebur's work again, we
generalize Duke's result to exact formulas for traces of singular
moduli of higher level modular functions.
\end{abstract}



\maketitle

\section{Introduction and statement of result}

 The classical $j$-invariant is defined for $z$ in the complex upper half
 plane $\mathcal{H}$ by
 $$j(z)= q^{-1} + 744 + 196884 q + \cdots, $$
 where $q=e(z)=e^{2 \pi i z}$ and $J(z)=j(z)-744$ is its normalized Hauptmodul
 for the group $\Gamma(1)=PSL_2(\mathbb Z)$.  All the modular
 groups discussed in this paper are subgroups of $\Gamma(1)$.
 For a positive integer $D$ congruent to 0 or 3 modulo 4, we
 denote by ${\mathcal Q}_D$ the set of positive definite integral binary quadratic
 forms $$Q(x,y)=[a,b,c]=ax^2+bxy+cy^2$$ with discriminant
 $-D=b^2-4ac$.  For each $Q\in {\mathcal Q}_D$, we let
 $$z_Q=\frac{-b+i\sqrt{D}}{2a},$$ the corresponding CM point in
 $\mathcal{H}$.  The group $\Gamma(1)$ acts on ${\mathcal Q}_D$ and we write
$\Gamma(1)_Q$ for the stabilizer of $Q$ in $\Gamma(1)$.
%That is, $\bar{\Gamma}(1)_Q=\{ \gamma\in
%\bar{\Gamma}(1)=PSL_2(\mathbb{Z})
  % \mid Q\circ \gamma=Q \}$.
The trace of $J$ singular moduli of discriminant $-D$ is defined
as
$$\textbf{t}_{J}(D)=\sum_{Q\in {\mathcal Q}_D/\Gamma(1)}
 \frac{1}{|{\Gamma}(1)_Q|} J(z_Q).$$
In \cite[Theorem 1]{Zagier}, Zagier proved the generating series
for the traces of singular moduli
 $$g(z):=q^{-1}-2-\mathop{\sum_{D>0}}_{D\equiv 0,3 (\textrm{mod} 4)}  \textbf{t}_{J}(D)q^D=q^{-1}-2+248q^3-492q^4+\cdots$$
 is a weakly holomorphic modular form of weight $3/2$ on
 $\Gamma_0(4)$.
% where the congruence subgroup $\Gamma_0(N)\subset\Gamma(1)$ is
% defined as $\Gamma_0(N)=\{ \sm a&b \\ c&d \esm \in SL_2(\mathbb Z)| c\equiv 0 \pmod{N}\}$.
Recently,  Bruinier, Jenkins, and Ono \cite{BJO}
 obtained an explicit formula for
 the Fourier coefficients of $g(z)$
 in terms of Kloosterman sums and Duke \cite{Duke} derived an exact formula for $\textbf{t}_{J}(D)$
 as follows;

 \begin{equation}\label{duke}
 \textbf{t}_{J}(D)=-24H(D)+\mathop{\sum_{c>0}}_{c\equiv 0
 (\textrm{mod} 4)}
 S_D(c)\sinh\left(\frac{4\pi\sqrt{D}}{c}\right),\end{equation}
 where $\displaystyle{H(D)=\sum_{Q\in {\mathcal Q}_D/\Gamma(1)}
 \frac{1}{|{\Gamma}(1)_Q|}}$ is the Hurwitz class number and

 $\displaystyle{S_D(c)=\sum_{x^2\equiv -D (\textrm{mod}\, c)} e(2x/c).}$
 Using these two results together,
 Duke \cite{Duke} reestablished Zagier's trace formula \cite[Theorem
 1]{Zagier}.

 The purpose of this paper is to give a generalization of (\ref{duke}) to
 traces of singular values of modular functions of any prime level
 $p$.  For prime $p$, let $\Gamma_0^*(p)$ be the
 group generated by
 $\Gamma_0(p)$ and the Fricke involution $W_p=\frac{1}{\sqrt{p}}\sm 0 & -1 \\  p&0 \esm$.
 Let $\mathcal{Q}_{D,p}$ denote the set of quadratic forms $Q\in \mathcal{Q}_D$ such that $a \equiv 0 \pmod p$
 on which $\Gamma_0^*(p)$ acts. Then the discriminant $-D$ is congruent to a square modulo $4p$.
 We choose an integer $\beta \pmod{2p}$ with $\beta^2 \equiv -D \pmod{4p}$ and consider
 the set $\mathcal{Q}_{D,p,\beta}
 =\{ [a,b,c]\in \mathcal{Q}_{D,p} \mid
 b \equiv \beta \pmod{2p} \}$ on which $\Gamma_0(p)$ acts.
 For a modular function $f$ for $\Gamma_0^*(p)$,
 we define the class number $H_p(D)$ (resp. $H_p^*(D)$) and
 the trace $\textbf{t}_{f}(D)$ (resp. $\textbf{t}_{f}^*(D)$) by
 \begin{eqnarray*}
  &H_p(D)=\underset{Q\in\mathcal{Q}_{D,p,\beta}/\Gamma_0(p)}{\sum}
  \frac{1}{|\Gamma_0(p)_Q|}; \quad
  & \textbf{t}_{f}(D) =
 \underset{Q\in\mathcal{Q}_{D,p,\beta}/\Gamma_0(p)}{\sum}
  \frac{1}{|\Gamma_0(p)_Q|} f (z_Q) \cr
  &H_p^*(D)=\underset{Q\in\mathcal{Q}_{D,p}/\Gamma_0^*(p)}{\sum}
  \frac{1}{|\Gamma_0^*(p)_Q|}; \quad
  & \textbf{t}_{f}^*(D) =
 \underset{Q\in\mathcal{Q}_{D,p}/\Gamma_0^*(p)}{\sum}
  \frac{1}{|\Gamma_0^*(p)_Q|} f (z_Q).
 \end{eqnarray*}
 Here $\Gamma_0(p)_Q$ and  $\Gamma^*_0(p)_Q$ are the stabilizers of
 $Q$ in $\Gamma_0(p)$ and $\Gamma^*_0(p)$, respectively.
 It is easy to see that
 \begin{equation} \label{class number}
 H_p^*(D) = \left\{%
\begin{array}{ll}
    \frac 12 H_p(D), & \hbox{ if $\beta\equiv 0$ or $p \pmod{2p},$} \\
    H_p(D), & \hbox{otherwise,} \\
\end{array}%
\right.
 \end{equation}
 and
 \begin{equation} \label{trace}
 \textbf{t}_{f}^*(D) = \left\{%
\begin{array}{ll}
    \frac 12 \textbf{t}_{f}(D), & \hbox{ if $\beta\equiv 0$ or $p \pmod{2p}$,} \\
    \textbf{t}_{f}(D), & \hbox{otherwise.} \\
\end{array}%
\right.
 \end{equation}
 The modularity for traces $\textbf{t}_{f}(D)$ was established
 by one of the authors in \cite{Kim1}, \cite{Kim2} when $\Gamma_0^*(p)$ is of genus zero.
 If $f$ is the Hauptmodul for such $\Gamma_0(p)^*$
 and if we define
 $\textbf{t}_f(-1)=-1$, $\textbf{t}_f(0)=2$ and
 $\textbf{t}_f(D)=0$
 for $D<-1$, then
 the series
 $\sum_{n,r} \textbf{t}_{f}(4pn-r^2) q^n \zeta^r$, where $\zeta=e(w)$ for a complex number $w$, is a
 weak Jacobi form of weight 2 and index $p$.
 Meanwhile, using theta correspondence, Bruinier and Funke
 \cite{BF} generalized Zagier's trace formula to traces of CM values of
 modular functions of arbitrary level. In particular, they showed that if $p$ is an odd
 prime and $f=\sum a(n)q^n$ is a modular function for
 $\Gamma_0^*(p)$ with $a(0)=0$, then
 $$ \mathop{\sum_{D>0}}_{-D\equiv \square \, (\textrm{mod} 4p)} \textbf{t}_f^*(D)q^D+\sum_{n\ge 1}(\sigma(n) + p\sigma(n/p))a(-n)
  -\sum_{m\ge 1}\sum_{n\ge 1} ma(-mn)q^{-m^2}$$
is a weakly holomorphic modular form of weight $3/2$ and level
$4p$.

We will obtain in the next section, the following exact formula
for $\textbf{t}_f^*(D)$ which is a generalization of (\ref{duke}).

\begin{thm}\label{exact formula}
Suppose $f$ is a modular function for $\Gamma_0^*(p)$ with
principal part $\sum_{m=1}^N a_m e(-mz)$ at $i\infty$ and define
for any positive integers $m$ and $c$,
$$S_D(m,c)=\sum_{x^2\equiv -D (\textrm{mod} \, c)} e(2mx/c).$$ Then
 $$
 \textbf{t}_f^*(D)=\sum_{m=1}^N a_m
 \left[
 c_m H_p^*(D)+\mathop{\sum_{c>0}}_{c\equiv 0\,(\textrm{mod}\, 4p)}
 S_D(m,c)\sinh\left(\frac{4\pi m\sqrt{D}}{c}\right)
 \right],
 $$
 where  $\displaystyle{c_m=-24\left( \frac{-p^{\alpha+1}}{p+1} \sigma(m/p^\alpha)
 +\sigma(m)\right)}$ with $p^\alpha|| m$.
\end{thm}

As an example, consider
 $$f=\left( \frac{\eta(z)}{\eta(37z)} \right)^2-2
 +37\left( \frac{\eta(37z)}{\eta(z)}\right)^2,$$
where $\eta(z)$ is the Dedekind eta function defined by
 $\eta(z)=q^{\frac1{24}}\prod_{n=1}^{\infty}(1-q^n)$.
 Then $f$ is a modular function for $\Gamma_0^*(37)$ which is of genus $1$ and has a
 Fourier expansion of the form
 $ q^{-3}-2q^{-2}-q^{-1}+0+O(q)$.
 Since the representatives for $\mathcal{Q}_{148,37,0}/\Gamma_0(37)$ are given by
 $[37,0,1]$ and $[74,-74,19]$, we find from equations (\ref{class number}),
 (\ref{trace}), and Theorem \ref{exact formula} that
 \begin{eqnarray*}
&&24\cdot \frac{3}{38}+\mathop{\sum_{c>0}}_{c\equiv 0
(\textrm{mod} 148)} \left[S_D(3,c)\sinh\left(\frac{12\pi
\sqrt{D}}{c}\right)
 -2 S_D(2,c)\sinh\left(\frac{8\pi \sqrt{D}}{c}\right)\right.\cr
 &&\quad \left. -S_D(1,c)\sinh\left(\frac{4\pi \sqrt{D}}{c}\right)\right]=\frac 12
 \left(f\left(\frac{\sqrt{37}i}{37}\right)+f\left(\frac{37+\sqrt{37}i}{74}\right)\right),
 \end{eqnarray*}
where the latter is known to be $-2$.

\section{Proof of Theorem \ref{exact formula}}

Throughout this section, $\Gamma$ denotes $\Gamma_0^*(p)$.
 For a positive integer $m$ we consider Niebur's Poincar\'e series
 \cite{Niebur}
 \begin{equation}\label{poincare}
 \mathcal F_m (z,s)=\sum_{M\in \Gamma_\infty\backslash\Gamma}
  e(-m {\rm Re}Mz) ({\rm Im} Mz)^{1/2} I_{s-1/2} (2\pi m {\rm Im}
  Mz),
 \end{equation}
 where $I_{s-1/2}$ is the modified Bessel function of the first kind.
 Then $\mathcal F_m (z,s)$ converges absolutely for ${\rm Re} \, s >
 1$ and satisfies
 \begin{equation} \label{harmonic}
  \mathcal F_m (Mz,s)=\mathcal F_m (z,s) \text{ for } M\in\Gamma
   \text{ and } \triangle \mathcal F_m (z,s)=s(1-s)\mathcal F_m
   (z,s),\end{equation}
   where $\triangle$ is the hyperbolic Laplacian
   $\triangle=-y^2(\partial^2_x+\partial^2_y)$ for $z=x+iy$.
 Niebur showed that $\mathcal F_m (z,s)$ has an analytic
 continuation to $s=1$ \cite[Theorem 5]{Niebur} and that
 $\mathcal F_m (z,s)$ has the following Fourier expansion
 \cite[Theorem 1]{Niebur}; for ${\rm Re} \, s > 1$,
 \begin{equation}\label{fourier}
 \mathcal F_m (z,s)=e(-mx)y^{1/2} I_{s-1/2}(2\pi m y)
  + \sum_{n=-\infty}^{\infty} b_n (y,s;-m) e(nx),
 \end{equation}
 where $b_n (y,s;-m) \to 0$ ($n\neq 0$) exponentially as $y\to
 i\infty$. Hence the pole of $\mathcal F_m (z,1)$ at $i\infty$
may occur only in
 $e(-mx)y^{1/2} I_{1/2}(2\pi m y)$, which is equal to
 \begin{equation} \label{pole}
 \frac{1}{\pi y^{1/2} m^{1/2}} \sinh (2\pi m y) y^{1/2} e(-mx)
 =\frac{1}{2\pi  m^{1/2}} (e(-mz)-e(-m\bar z)).
 \end{equation}
 So if we multiply $\mathcal F_m (z,1)$ by $2\pi  m^{1/2}$,
 then the coefficient of $e(-mz)$ is normalized. Now we need to
 compute the constant term in $(2\pi  m^{1/2})\mathcal F_m (z,1)$.
 It follows from \cite[Theorem 1]{Niebur} that
 $b_0(y,s,-m)=a_m(s) y^{1-s}/(2s-1)$. Here
 \begin{equation}\label{def}
 a_m(s)=2\pi^s m^{s-1/2} \phi_m(s)/\Gamma(s) \text{ \quad and \quad }
 \phi_m(s)=\sum_{c>0} S(m,0;c) c^{-2s},
 \end{equation}
 where $S(m,n;c)$ is the general Kloostermann sum
 $\sum_{0\le d < |c|} e((ma+nd)/c)$ for $\sm a & * \\ c&d \esm\in
 \Gamma$.
 Note that if $M=\sm a&b\\ c&d \esm\in \Gamma=\Gamma_0^*(p)$, then
 $M\in\Gamma_0(p)$ or $M$ is of the form
 $\sm \sqrt{p}x & y/\sqrt p \\ \sqrt p z& \sqrt p w \esm_{\det=1}$
 with $x,y,z,w\in \mathbb Z$.
 In the former case, $c$ is a multiple of $p$ and in the latter
 case, $c=\sqrt{p}z$ with $p\nmid z$. For $n\in \mathbb Z^+$,
 let $u_m(n)$ denote the sum of $m$-th powers of primitive $n$-th
 roots of unity. We observe that
 $$S(m,0;c)=\left\{%
\begin{array}{ll}
    u_m(c), & \hbox{ if $p\mid c$,} \\
    u_m(z), & \hbox{ if $c=\sqrt{p}z$ with $p\nmid z$.} \\
\end{array}%
\right.$$
If we define $u_m^*(n)=\left\{%
\begin{array}{ll}
    u_m(n), & \hbox{ if $p\mid n$,} \\
    p^{-s} u_m(n), & \hbox{ if $p\nmid n$.} \\
\end{array}%
\right.,$ then
 \begin{eqnarray} \label{phi_m}
 &&p^s \phi_m(s) \zeta (2s)
  =p^s \sum_{c>0} S(m,0;c) c^{-2s} \sum_{c'\in \mathbb{Z}^+}
  c'^{-2s}\cr
 &&\hspace{1cm} =\sum_{c\in \mathbb{Z}^+} (p^s u_m^*(c)) c^{-2s}
   \sum_{c'\in \mathbb{Z}^+} c'^{-2s}
  = \sum_{k\in \mathbb{Z}^+} \left(\sum_{c\mid k} p^s u_m^*(c)\right) k^{-2s}.
 \end{eqnarray}
Note that if $p\nmid k$, then
 \begin{equation} \label{nmid}
 \sum_{c\mid k} p^s u_m^*(c)=\sum_{c\mid k} u_m(c)=
 \left\{%
\begin{array}{ll}
    k, & \hbox{ if $k\mid m$,} \\
    0, & \hbox{ if $k\nmid m$.} \\
\end{array}%
\right.
 \end{equation}
and if $k=p^l k'$ with $l\ge 1$ and $p\nmid k'$, then
 \begin{equation}\label{pdivk}
 \sum_{c\mid k} p^s u_m^*(c)
 =\sum_{d\mid k'} p^s u_m^*(d) +
 \mathop{\sum_{c\mid k}}_{p\mid c} p^s u_m^*(c)
 =\sum_{d\mid k'}  u_m(d) +
\mathop{\sum_{c\mid k}}_{p\mid c} p^s u_m(c).
 \end{equation}
 By adding $(p^s-1)\sum_{d\mid k'}  u_m(d)$ on both sides of (\ref{pdivk}), we obtain
 \begin{equation*}
 (p^s-1)\sum_{d\mid k'}  u_m(d)+ \sum_{c\mid k} p^s u_m^*(c)
 = \sum_{c\mid k} p^s u_m(c).
 \end{equation*}
 Since
 $\displaystyle{\sum_{d\mid k'} u_m(d)=
 \left\{%
\begin{array}{ll}
    k', & \hbox{ if $k'\mid m$,} \\
    0, & \hbox{ if $k'\nmid m$} \\
\end{array}%
\right. }$\quad  and \quad
 $\displaystyle{\sum_{c\mid k} p^s
u_m(c)=
 \left\{%
\begin{array}{ll}
    p^s k, & \hbox{ if $k\mid m$,} \\
    0, & \hbox{ if $k\nmid m$} \\
\end{array}%
\right. }$,

we find that
 \begin{equation} \label{mid}
 \sum_{c\mid k} p^s u_m^*(c) =
 \left\{%
\begin{array}{ll}
    p^s k+(1-p^s)k', & \hbox{ if $k\mid m$,} \\
    (1-p^s)k', & \hbox{ if $k\nmid m$ and $k'\mid m$,} \\
    0, & \hbox{ if $k\nmid m$ and $k'\nmid m$.} \\
\end{array}%
\right.
 \end{equation}
Writing $m=p^\alpha m'$ with $p\nmid m'$, we can deduce from
(\ref{nmid}) and (\ref{mid}) that
 \begin{eqnarray}\label{general}
 &&\sum_{k\in\mathbb Z_+} (\sum_{c\mid k} p^s u_m^*(c)) k^{-2s}
 =\sum_{k' \mid m'} k' k'^{-2s}
   +\sum_{l=1}^\infty \sum_{k' \mid m'} (1-p^s)k' (p^l
   k')^{-2s}\cr
&& \hspace{5cm}  +\sum_{l=1}^\alpha \sum_{k' \mid m'} p^s(p^l k')
(p^l k')^{-2s}\cr
 &&= \sigma_{1-2s}(m')+(1-p^s)\sigma_{1-2s} (m')\sum_{l=1}^\infty
 (p^{-2s})^l +p^s \mathop{\sum_{1\le l \le \alpha}}{k'\mid m'} (p^l
 k')^{1-2s} \cr
 &&=\sigma_{1-2s}(m')\left[ 1+(1-p^s)\frac{p^{-2s}}{1-p^{-2s}}\right]+p^s
 (\sigma_{1-2s}(m)-\sigma_{1-2s}(m')) \cr
 &&=\frac{-p^{2s}}{1+p^s} \sigma_{1-2s}(m/p^\alpha)+p^s
 \sigma_{1-2s}(m)
 \end{eqnarray}
Recall the constant term in $(2\pi  m^{1/2})\mathcal F_m (z,1)$ is
$$\lim_{s\to 1} 2\pi  m^{1/2} b_0(y,s,-m)=
 \lim_{s\to 1} 2\pi  m^{1/2} a_m(s) y^{1-s}/(2s-1).$$
By the definition of $a_m(s)$ in (\ref{def}), it is equal to
 $$\lim_{s\to 1} 2\pi m^{1/2} (2\pi^s m^{s-1/2}
 \phi_m(s)/\Gamma(s)) y^{1-s}/(2s-1).$$
It follows from (\ref{phi_m}) and (\ref{general}) that this limit
goes to $$\frac{4\pi^2 m}{p\zeta(2)} \left( \frac{-p^{2}}{1+p}
\sigma_{-1}(m/p^\alpha)+p
 \sigma_{-1}(m)\right).$$
%$$ 24 \left( \frac{-p}{1+p}m\sigma_{-1}(m/p^\alpha)+m\sigma_{-1}(m)\right)$$
Thus simple calculations lead us to have the constant term of
$(2\pi  m^{1/2})\mathcal F_m (z,1)$ as
\begin{equation}\label{constant}
 24 \left(\frac{-p^{\alpha+1}}{1+p}
  \sigma(m/p^\alpha)+\sigma(m)\right)=-c_m. \end{equation}

 Now we define $$\mathcal F_m^*(z,s)=(2\pi  m^{1/2})\mathcal F_m
 (z,s)+c_m.$$ Then by (\ref{harmonic}), (\ref{fourier}), and (\ref{pole}),
 $\mathcal F_m^*(z,1)$ is $\Gamma$-invariant harmonic function and
 $\mathcal F_m^*(z,1)-e(-mz)$ has a zero at $i\infty$.
 Hence it follows from \cite[Theorem 6]{Niebur} that
 \begin{equation*}
 f(z)=\sum_{m=1}^N a_m \mathcal F_m^*(z,1)
 \end{equation*}
 for any modular function $f$ for
 $\Gamma_0^*(p)$ with principal part $\sum_{m=1}^N a_m e(-mz)$ at
 $i\infty$.  Hence
 \begin{equation} \label{Niebur}\textbf{t}_f^*(D)=\sum_{m=1}^N a_m
 \left(\sum_{Q\in\mathcal Q_{D,p}/\Gamma} \frac{1}{|\Gamma_{Q}|}
 \mathcal F_m^*(z_{Q},1)\right).
 \end{equation}  In order to determine the value
 $\displaystyle{\sum_{Q\in\mathcal Q_{D,p}/\Gamma} \frac{1}{|\Gamma_{Q}|}
 \mathcal F_m^*(z_{Q},1)}$, we first compute for ${\rm Re}\, s>1$,
 \begin{equation}\label{F*}\sum_{Q\in\mathcal Q_{D,p}/\Gamma} \frac{1}{|\Gamma_{Q}|}
 \mathcal F_m^*(z_{Q},s)={c_m H_p^*(D)
   +2\pi\sqrt{m}\sum_{Q\in\mathcal Q_{D,p}/\Gamma} \frac{1}{|\Gamma_{Q}|}
 \mathcal F_m(z_{Q},s).}\end{equation}
By Poincar\'{e} series expansion of $\mathcal F_m(z_{Q},s)$ in
(\ref{poincare}),
$$2\pi\sqrt{m}\sum_{Q\in\mathcal Q_{D,p}/\Gamma}\frac{ \mathcal F_m(z_{Q},s)}{|\Gamma_{Q}|}
=2\pi\sqrt{m}
 \sum_{Q\in\mathcal Q_{D,p}/\Gamma_\infty}
 e(-m {\rm Re}\,z_Q) ({\rm Im} z_Q)^{1/2} I_{s-1/2} (2\pi m {\rm Im}
  z_Q).$$
The series in the latter is equal to
\begin{eqnarray*}
&& \sum_{[ap,b,c]\in\mathcal Q_{D,p}/\Gamma_\infty}
 e\left(\frac{2mb}{4pa}\right) \left(\frac{2\sqrt D}{4pa}\right)^{1/2}
 I_{s-1/2} \left(2\pi m \frac{2\sqrt D}{4pa}\right) \cr
&&= \sum_{a=1}^\infty\mathop{\sum_{x \, (\textrm{mod} \,
2ap)}}_{x^2\equiv -d \, (\textrm{mod}\,  4ap)}
 e\left(\frac{2mx}{4pa}\right) \left(\frac{2\sqrt D}{4pa}\right)^{1/2}
 I_{s-1/2} \left(2\pi m \frac{2\sqrt D}{4pa}\right) \cr
&&= \mathop{\sum_{c>0}}_{c\equiv 0 \, (\textrm{mod}\, 4p)}
 \frac 12 S_D(m,c) \left(\frac{2\sqrt D}{c}\right)^{1/2}
 I_{s-1/2} \left(2\pi m \frac{2\sqrt D}{c}\right),
 \end{eqnarray*}
 which converges uniformly for $s\in [1,2]$. Therefore, by
 (\ref{F*}),
$$\sum_{Q\in\mathcal Q_{D,p}/\Gamma} \frac{1}{|\Gamma_{Q}|}
 \mathcal F_m^*(z_{Q},1)= c_m H_p^*(D)+
 \mathop{\sum_{c>0}}_{c\equiv 0 \, (\textrm{mod} \, 4p)} S_D(m,c)\sinh\left(\frac{4\pi m \sqrt
 D}{c}\right).$$
 This combined with (\ref{Niebur}) completes the proof.

\end{document}